\documentclass[a4paper,12pt]{article}
\usepackage{amscd}
\usepackage{amsmath,amsfonts,amssymb,amscd}
\usepackage{indentfirst,graphicx,epsfig}
\usepackage{graphicx,psfrag}
\usepackage{ifpdf}
\usepackage{float}

\input{epsf}

\newtheorem{thm}{Theorem}

\newtheorem{lem}{Lemma}

\def\qed{\hfill \nopagebreak\rule{5pt}{8pt}}
\def\pf{\noindent {\it Proof.} }

\title {\bf The complexity of determining the\\ rainbow vertex-connection
of graphs\footnote{Supported by NSFC and ``the Fundamental Research
Funds for the Central Universities". }}
\author{
\small Lily Chen, Xueliang Li, Yongtang Shi\\
\small Center for Combinatorics and LPMC-TJKLC \\
\small Nankai University, Tianjin 300071, China \\
\small Email: lily60612@126.com, lxl@nankai.edu.cn,
shi@nankai.edu.cn
\date{}}
\begin{document}
\maketitle
\begin{abstract}
A vertex-colored graph is {\it rainbow vertex-connected} if any two
vertices are connected by a path whose internal vertices have
distinct colors, which was introduced by Krivelevich and Yuster. The
{\it rainbow vertex-connection} of a connected graph $G$, denoted by
$rvc(G)$, is the smallest number of colors that are needed in order
to make $G$ rainbow vertex-connected. In this paper, we study the
computational complexity of vertex-rainbow connection of graphs and
prove that computing $rvc(G)$ is NP-Hard. Moreover, we show that it
is already NP-Complete to decide whether $rvc(G)=2$. We also prove
that the following problem is NP-Complete: given a vertex-colored
graph $G$, check whether the given coloring makes $G$ rainbow
vertex-connected.\\
[2mm] Keywords: coloring; rainbow vertex-connection; computational complexity\\
[2mm] AMS Subject Classification (2010): 05C15, 05C40, 68Q25, 68R10.
\end{abstract}

\section{Introduction}

All graphs considered in this paper are simple, finite and
undirected. We follow the notation and terminology of Bondy and
Murty \cite{BM}.

An edge-colored graph is {\it rainbow connected} if any two vertices
are connected by a path whose edges have distinct colors. This
concept of rainbow connection in graphs was introduced by Chartrand
et al. in \cite{CJMZ}. The {\it rainbow connection number} of a
connected graph $G$, denoted by $rc(G)$, is the smallest number of
colors that are needed in order to make $G$ rainbow connected.
Observe that $diam(G)\leq rc(G)\leq n-1$, where $diam(G)$ denotes
the diameter of $G$. It is easy to verify that $rc(G)=1$ if and only
if $G$ is a complete graph, that $rc(G)=n-1$ if and only if $G$ is a
tree. There are some approaches to study the bounds of $rc(G)$, we
refer to \cite{CLRTY,KY,S}.

In \cite{KY}, Krivelevich and Yuster proposed the concept of rainbow
vertex-connection.  A vertex-colored graph is {\it rainbow
vertex-connected} if any two vertices are connected by a path whose
internal vertices have distinct colors. The {\it rainbow
vertex-connection} of a connected graph $G$, denoted by $rvc(G)$, is
the smallest number of colors that are needed in order to make $G$
rainbow vertex-connected. An easy observation is that if $G$ is of
order $n$ then $rvc(G)\leq n-2$ and $rvc(G)=0$ if and only if $G$ is
a complete graph. Notice that $rvc(G)\geq diam(G)-1$ with equality
if the diameter is $1$ or $2$. For rainbow connection and rainbow
vertex-connection, some examples are given to show that there is no
upper bound for one of parameters in terms of the other in
\cite{KY}. Krivelevich and Yuster \cite{KY} proved that if $G$ is a
graph with $n$ vertices and minimum degree $\delta$, then
$rvc(G)<11n/\delta$. In \cite{LS}, the authors improved this bound
for given order $n$ and minimum degree $\delta$.

Besides its theoretical interest as being a natural combinatorial
concept, rainbow connectivity can also find applications in
networking. Suppose we want to route messages in a cellular network
such that each link on the route between two vertices is assigned
with a distinct channel. The minimum number of used channels is
exactly the rainbow connection of the underlying graph.

The computational complexity of rainbow connection has been studied.
In \cite{CLRTY}, Caro et al. conjectured that computing $rc(G)$ is
an NP-Hard problem, as well as that even deciding whether a graph
has $rc(G)=2$ is NP-Complete. In \cite{CFMY}, Chakraborty et al.
confirmed this conjecture. Motivated by the proof of \cite{CFMY}, we
consider the computational complexity of rainbow vertex-connection
$rvc(G)$ of graphs. It is not hard to image that this problem is
also NP-hard, but a rigorous proof is necessary. This paper is to
give such a proof, which follows a similar idea of \cite{CFMY}, but
by reducing $3$-SAT problem to some other new problems, that
computing $rvc(G)$ is NP-Hard. Moreover, we show that it is already
NP-Complete to decide whether $rvc(G)=2$. We also prove that the
following problem is NP-Complete: given a vertex-colored graph $G$,
check whether the given coloring makes $G$ rainbow vertex-connected.

\section{Rainbow vertex-connection.}

For two problems $A$ and $B$, we write $A \preceq B$, if problem $A$
is polynomially reducible to problem $B$. Now, we give our first
theorem.
\begin{thm}\label{theorem2}
The following problem is NP-Complete: given a vertex-colored graph
$G$, check whether the given coloring makes $G$ rainbow
vertex-connected.
\end{thm}

Now we define Problem 1 and Problem 2 as follows. We will prove
Theorem \ref{theorem2} by reducing Problem 1 to Problem 2, and then
$3$-SAT problem to Problem 1.

{\bf Problem 1}~ $s-t$ rainbow vertex-connection.

\noindent Given: Vertex-colored graph $G$ with two vertices $s,\ t$.\\
Decide: Whether there is a rainbow vertex-connected path between $s$
and $t$?

{\bf Problem 2}~ Rainbow vertex-connection.

\noindent Given: Vertex-colored graph $G$.\\
Decide: Whether $G$ is rainbow vertex-connected under the coloring?

\begin{lem}
Problem 1 $\preceq$ Problem 2.
\end{lem}

\pf Given a vertex-colored graph $G$ with two vertices $s$ and $t$.
We want to construct a new graph $G'$ with a vertex coloring such
that $G'$ is rainbow vertex-connected if and only if there is a
rainbow vertex-connected path from $s$ to $t$ in $G$.

Let $V=\{v_1, v_2,\ldots,v_{n-1},v_n\}$ be vertices of $G$, where
$v_1=s$ and $v_n=t$. We construct $G'$ as follows. Set
$$V'=V\cup\{s',t',a,b\}$$ and
$$E'=E\cup \{s's,t't\}\cup\{av_i,bv_i:i\in [n]\}.$$

Let $c$ be the vertex coloring of $G$, we define the vertex coloring
$c'$ of $G'$ by $c'(v_i)=c(v_i)$ for $i\in \{2,3,\ldots,n-1\}$,
$c'(s)=c'(a)=c_1,$ $c'(t)=c'(b)=c_2$, where $c_1,c_2$ are the two
new colors.

Suppose $c'$ makes $G'$ rainbow vertex-connected. Since each path
$Q$ from $s'$ to $t'$ must go through $s$ and $t$, $Q$ can not
contain $a$ and $b$ as $c'(s)=c'(a)=c_1$ and $c'(t)=c'(b)=c_2$.
Therefore, any rainbow vertex-connected path from $s'$ to $t'$ must
contain a rainbow vertex-connected path from $s$ to $t$ in $G$. Thus
there is a rainbow vertex-connected path from $s$ to $t$ in $G$
under the coloring $c$.

Now assume that there is a rainbow vertex-connected path from $s$ to
$t$ in $G$ under the coloring $c$. To prove that $G'$ is rainbow
vertex-connected. First, the rainbow vertex-connected path from $s'$
to $v_i$ can be formed by going through $s$ and $b$, then to $v_i$
for $i\in\{2,3,\ldots,n\}$. The rainbow vertex-connected path from
$s'$ to $t'$ can go through $s$ and $t$ and a rainbow
vertex-connected path from $s$ to $t$ in $G$. The rainbow
vertex-connected path from $t'$ to $v_i$ can be formed by going
through $t$ and $a$, then to $v_i$ for $i\in\{2,3,\ldots,n\}$. For
the other pairs of vertices, there is a path between them with
length less than $3$, thus they are obvious rainbow
vertex-connected.\qed

\begin{lem}
$3$-SAT $\preceq$ Problem 1.
\end{lem}

\pf Let $\phi$ be an instance of $3$-SAT with clauses
$c_1,c_2,\ldots, c_m$ and variables $x_1,x_2,\ldots,x_n$. We
construct a graph $G_\phi$ with special vertices $s$ and $t$.

First, we introduce $k$ new vertices $v_1^j,v_2^j,\ldots,v_k^j$ for
each $x_j\in c_i$ and $\ell$ new vertices
$\overline{v}_1^j,\overline{v}_2^j,\ldots,\overline{v}_\ell^j$ for
each $\overline{x}_j\in c_i$. Without loss of generality, we assume
that $k\geq1$ and $\ell\geq1$, otherwise $\phi$ can be simplified.

Next, for each $v_a^j,\ a\in [k]$, we introduce $\ell$ new vertices
$v_{a1}^j, v_{a2}^j,\ldots,v_{a\ell}^j$, which form a path in this
order. Similarly, for each $\overline{v}_b^j,\ b\in [\ell]$, we
introduce $k$ new vertices $\overline{v}_{1b}^j,
\overline{v}_{2b}^j,\ldots,\overline{v}_{kb}^j$, which form a path
in that order. Therefore, for $x_j\in c_i$, there are $k$ paths of
length $\ell-1$, and for $\overline{x}_j\in c_i$, there are $\ell$
paths of length $k-1$. For each path $Q$ in $c_i$ ($i\in [m]$), we
join the original vertex of $Q$ to the terminal vertices of all
paths in $c_{i-1}$, where $c_0$ is the vertex $s$. And for each path
in $c_m$, we join its terminal vertex to $t$. Thus, a new graph
$G_\phi$ is obtained.

Now we define a vertex coloring of $G_\phi$. For every variable
$x_j$, we introduce $k\times \ell$ distinct colors
$\alpha_{1,1}^j,\alpha_{1,2}^j,\ldots,\alpha_{k,\ell}^j$\,. We color
vertices $v_{a1}^j, v_{a2}^j,\ldots,v_{a\ell}^j$ with colors
$\alpha_{a,1}^j,\alpha_{a,2}^j,\ldots,\alpha_{a,\ell}^j$\,,
respectively, and color $\overline{v}_{1b}^j,
\overline{v}_{2b}^j,\ldots,\overline{v}_{kb}^j$ with colors
$\alpha_{1,b}^j,\alpha_{2,b}^j,\ldots,\alpha_{k,b}^j$\,,
respectively, where $a\in[k]$ and $b\in[\ell]$.

If $G_\phi$ contains a rainbow vertex-connected $s-t$ path $Q$, then
$Q$ must contain one of the newly built paths in each $c_i$,
$i\in[m]$, and the path $v_{a1}^jv_{a2}^j\ldots v_{a\ell}^j$ and
$\overline{v}_{1b}^j\overline{v}_{2b}^j\ldots \overline{v}_{kb}^j$
can not both appear in $Q$. If $v_{a1}^jv_{a2}^j\ldots v_{a\ell}^j$
appears in $Q$, set $x_j=1$, and if
$\overline{v}_{1b}^j\overline{v}_{2b}^j\ldots \overline{v}_{kb}^j$
appears in $Q$, set $x_j=0$. Clearly, we have $\phi=1$ in this
assignment.\qed

\section{Rainbow vertex-connection $2$.}

Before proceeding, we first define three problems.

{\bf Problem 3}~ Rainbow vertex-connection $2$.

\noindent Given: Graph $G=(V,E)$.\\ Decide: Whether there is a
vertex coloring of $G$ with two colors such that all pairs $(u,v)\in
V(G)\times V(G)$ are rainbow vertex-connected?

{\bf Problem 4}~ Subset rainbow vertex-connection $2$.

\noindent Given: Graph $G=(V,E)$ and a set of pairs $P\subseteq
V(G)\times V(G).$\\
Decide: Whether there is a vertex coloring of $G$ with two colors
such that all pairs $(u,v)\in P$ are rainbow vertex-connected?

{\bf Problem 5}~ Different subsets rainbow vertex-connection 2.

\noindent Given: Graph $G=(V,E)$ and two disjoint subsets $V_1,V_2$
of $V$ with a one to one corresponding $f:V_1\rightarrow V_2$.\\
Decide: Whether there is a vertex coloring of $G$ with two colors
such that $G$ is rainbow vertex-connected and for each $v\in V_1$,
$v$ and $f(v)$ are assigned different colors.

In the following, we will reduce Problem 4 to Problem 3 and then
reduce Problem 5 to Problem 4. Finally, we will show it is
NP-Complete to decide whether $rvc(G)=2$ by reducing $3$-SAT problem
to Problem 3.

\begin{lem}
Problem 4  $\preceq$ Problem 3.
\end{lem}

\pf Given a graph $G=(V,E)$ and a set of pairs $P\subseteq
V(G)\times V(G)$, we construct a graph $G'=(V',E')$ as follows.

For each vertex $v\in V$, we introduce a new vertex $x_v$; for every
pair $(u,v)\in (V\times V)\setminus P$, we introduce two new
vertices $x_{(u,v)}^1$ and $x_{(u,v)}^2$; we also add two new
vertices $s,t$. Set
$$V'=V\cup\{x_v: v\in V\}\cup \{x_{(u,v)}^1, x_{(u,v)}^2:
(u,v)\in(V\times V)\setminus P\} \cup \{s,t\}$$ and\\
$E'=E\cup\{vx_v: v\in V\}\cup
\{ux_{(u,v)}^1,x_{(u,v)}^1x_{(u,v)}^2,x_{(u,v)}^2v:(u,v)\in(V\times
V)\setminus P \} \cup \{sx_{(u,v)}^1,tx_{(u,v)}^2:(u,v)\in(V\times
V)\setminus P\}\cup\{sx_v,tx_v:v\in V\}.$

\noindent Observe that $G$ is a subgraph of $G'$. In the following,
we will prove that $G'$ is $2$-rainbow vertex-connected if and only
if there is a vertex coloring of $G$ with two colors such that all
pairs $(u,v)\in P$ are rainbow vertex-connected.

Now suppose there is a vertex coloring of $G'$ with two colors which
makes $G'$ rainbow vertex-connected. For each pair $(u,v)\in P$, the
paths of length no more than $3$ that connects $u$ and $v$ have to
be in $G$. Thus, with the coloring all pairs in $P$ are rainbow
vertex-connected. On the other hand, let $c: V\rightarrow \{1,2\}$
be one coloring of $G$ such that all pairs $(u,v)\in P$ are rainbow
vertex-connected. We extend the coloring as follows: $c(x_v)=1$ for
all $v\in P$, $c(x_{(u,v)}^1)=1$ and $c(x_{(u,v)}^2)=2$ for all
$(u,v)\in(V\times V)\setminus P$, $c(s)=c(t)=2$. We can see that
$G'$ is indeed rainbow vertex-connected under this coloring. \qed

\begin{lem}
Problem 5 $\preceq$ Problem 4.
\end{lem}

\pf Given a graph $G=(V,E)$ and two disjoint subsets $V_1,V_2$ of
$V$ with a one to one corresponding $f$. Assume that
$V_1=\{v_1,v_2,\ldots, v_k\}$ and $V_2=\{w_1,w_2,\ldots,w_k\}$
satisfying that $w_i=f(v_i)$ for each $i\in[k]$. We construct a new
graph $G'=(V',E')$ as follows.

We introduce six new vertices
$x_{v_iw_i}^1,x_{v_iw_i}^2,x_{v_iw_i}^3,x_{v_iw_i}^4,
x_{v_iw_i}^5,x_{v_iw_i}^6$ for each pair $(v_i,w_i)$, $i\in [n]$. We
add a new vertex $s$. Set
$$V'=V\cup \{x_{v_iw_i}^j:i\in[k],j\in [6]\}\cup {\{s\}},$$
and\\
$E'=E\cup\{sx_{v_iw_i}^5,x_{v_iw_i}^5v_i,v_ix_{v_iw_i}^1,
x_{v_iw_i}^1x_{v_iw_i}^2,x_{v_iw_i}^2x_{v_iw_i}^3,x_{v_iw_i}^3x_{v_iw_i}^4,
x_{v_iw_i}^4w_i$, $w_ix_{v_iw_i}^6,x_{v_iw_i}^6s:i\in[k]\}$.

\noindent We define $P$ by:\\
$P=\{(u,v):u,v\in V\}\cup\{(x_{v_iw_i}^5,x_{v_iw_i}^2),
(v_i,x_{v_iw_i}^3),(x_{v_iw_i}^1,x_{v_iw_i}^4),(x_{v_iw_i}^2,w_i),
(x_{v_iw_i}^3,x_{v_iw_i}^6):i\in [k]\}$.

Suppose there is a vertex coloring of $G'$ with two colors such that
all pairs $(u,v)\in P$ are rainbow vertex-connected. Observe that
$G$ is a subgraph of $G'$. For all $(u,v)\in V\times V$, they are
belong to $P$ and the paths connect them with length no more than
$3$ are belong to $G$, thus $G$ is rainbow vertex-connected.  We
have $c(v_i)\neq c(w_i)$, since $\{(x_{v_iw_i}^5,x_{v_iw_i}^2),
(v_i,x_{v_iw_i}^3),(x_{v_iw_i}^1,x_{v_iw_i}^4),(x_{v_iw_i}^2,w_i),
(x_{v_iw_i}^3,x_{v_iw_i}^6):i\in [k]\}$ are rainbow vertex-connected
in $G'$.

On the other hand, if there is a $2$-vertex coloring $c$ of $G$ such
that $G$ is rainbow vertex-connected and $v_i, w_i$ are colored
differently, we color $G'$ with coloring $c'$ as follows. For $v\in
V$, $c'(v)=c(v)$. If $c(v_i)=1,\ c(w_i)=2$, then
$c'(x_{v_iw_i}^1)=c'(x_{v_iw_i}^3)=2,$
$c'(x_{v_iw_i}^2)=c'(x_{v_iw_i}^4)=1.$ If $c(v_i)=2,\ c(w_i)=1$,
then $c'(x_{v_iw_i}^1)=c'(x_{v_iw_i}^3)=1,$
$c'(x_{v_iw_i}^2)=c'(x_{v_iw_i}^4)=2.$ For all other vertices, we
assign them by color $1$ or $2$ arbitrarily. It is easy to check
that all $(u,v)\in P$ are rainbow vertex-connected. \qed

\begin{lem}
$3$-SAT $\preceq$ Problem 5.
\end{lem}

\pf Let $\phi$ be an instance of $3$-SAT with clauses
$c_1,c_2,\ldots, c_m$ and variables $x_1,x_2,\ldots,x_n$. We
construct a new graph $G_\phi$ and define two disjoint vertex sets
with a one to one corresponding $f$. Add two new vertices $s,t$. Set\\
$$V_\phi=\{c_i:i\in [m]\}\cup \{x_i,\overline{x}_i:i\in
[n]\}\cup\{s,t\}$$ and\\
$E_\phi=\{c_ic_j:i,j \in [m]\}\cup\{tx_i,t\overline{x}_i:i\in
[n]\}\cup \{x_ic_j:x_i\in
c_j\}\cup\{\overline{x}_ic_j:\overline{x}_i\in c_j\}\cup \{st\}$.

\noindent We define $V_1=\{x_1,x_2,\ldots,x_n\},$
$V_2=\{\overline{x}_1,\overline{x}_2,\ldots,\overline{x}_n\}$ and
$f:V_1\rightarrow V_2$ satisfying that $f(x_i)=\overline{x}_i$. Now
we show that $G_\phi$ is 2-rainbow vertex-connected with different
colors between $x_i$ and $\overline{x}_i$ if and only if $\phi$ is
satisfiable.

Suppose there is a vertex coloring $c:V_\phi\rightarrow \{0,1\}$
such that $G_\phi$ is rainbow vertex-connected and $x_i$,
$\overline{x}_i$ are colored differently. We first suppose $c(t)=0$,
set the value of $x_i$ as the corresponding color of $x_i$. For each
$i$, consider the rainbow vertex-connected path $Q$ between vertices
$s$ and $c_i$, there must exist some $j$ such that we can write
$Q=stx_jc_i$ or $Q=st\overline{x}_jc_i$. Without loss of generality,
suppose $Q=stx_jc_i$. Since $c(t)=0$, we have $c(x_j)=1$. Thus, the
value of $x_j$ is $1$, which concludes that $c_i=1$ as $x_j\in c_i$
by the construction of $G_\phi$. For the other case $c(t)=1$, we set
$x_i=1$ if $c(x_i)=0$ and $x_i=0$ otherwise. By some similar
discussions, we also have $\phi=1$.

On the other hand, for a given truth assignment of $\phi$, we color
$G_\phi$ as follows: $c(t)=0$ and $c(c_i)=1$ for $i\in[m]$; if
$x_i=1$, then $c(x_i)=1$ and $c(\overline{x}_i)=0$; otherwise,
$c(x_i)=0$ and $c(\overline{x}_i)=1$. We can easily check that
$G_\phi$ is rainbow vertex-connected.\qed

From the above three lemmas, we conclude our second theorem.
\begin{thm}
Given a graph $G$, deciding whether $rvc(G)=2$ is NP-Complete. Thus,
computing $rvc(G)$ is NP-Hard.
\end{thm}

\end{document}